\documentclass[aos, reqno, preprint, 11pt]{article}%

\usepackage{pdfsync}

\RequirePackage{amsthm, amsmath, amsfonts, amssymb}%

\usepackage{graphicx, color}%
\usepackage{tikz}%
\usepackage{mathrsfs}

\usepackage{tikz}

\definecolor{darkblue}{rgb}{0.0,0.0,0.7}

\newcommand{\norm}[1]{\|#1\|}%
\newcommand{\ind}[1]{\mathbf 1{(#1)}}%

\renewcommand{\P}{\mathbb P}

\newcommand{\cH}{{\mathcal H}}
\newcommand{\G}{{\mathbf G}}
\renewcommand{\H}{{\mathbf H}}

\newcommand{\vt}{\vartheta}

\DeclareMathOperator*{\argmin}{argmin}

\DeclareMathOperator{\pen}{pen}

\DeclareMathOperator*{\re}{RE}

\newcommand{\1}{{\rm 1}\kern-0.24em{\rm I}}

\newtheorem{theorem}{Theorem}%
\newtheorem{corollary}{Corollary}%
\newtheorem{lemma}{Lemma}%
\theoremstyle{assumption}%
\newtheorem{assumption}{Assumption}{\bf}{\rm}%
\theoremstyle{remark}%
\newtheorem{remark}{Remark}%

\newcommand{\eps}{\epsilon}

\newcommand{\bs}{\boldsymbol}

\newcommand{\R}{\mathbb R}

\newcommand{\inr}[1]{\langle #1 \rangle}

\newcommand{\grad}{\nabla}
\newcommand{\cp}{\complement}
\newcommand{\sgn}{{\rm sgn}\,}
\newcommand{\C}{\mathbb C}

\newcommand{\RE}{{\rm RE}}

\begin{document}

\title{High-dimensional additive hazards models and the Lasso}

\author{St\'ephane Ga\"iffas$^{1}$ \and Agathe Guilloux$^{1,2}$}

\footnotetext[1]{Universit\'e Pierre et Marie Curie - Paris~6,
  Laboratoire de Statistique Th\'eorique et Appliqu\'ee. \emph{email}:
  \texttt{stephane.gaiffas@upmc.fr}}

\footnotetext[2]{Centre de Recherche Saint-Antoine (UMR S 893)
  \emph{email}: \texttt{agathe.guilloux@upmc.fr}}

\footnotetext[3]{This work is supported by French Agence Nationale de la
  Recherche (ANR) ANR Grant \textsc{``Prognostic''}
  ANR-09-JCJC-0101-01.}

\maketitle

\begin{abstract}
  We consider a general high-dimensional additive hazards model in a
  non-asymptotic setting, including regression for censored-data. In
  this context, we consider a Lasso estimator with a fully data-driven
  $\ell_1$ penalization, which is tuned for the estimation problem at
  hand. We prove sharp oracle inequalities for this estimator. Our
  analysis involves a new ``data-driven'' Bernstein's inequality, that
  is of independent interest, where the
  predictable variation is replaced by the optional variation. \\

  \noindent%
\end{abstract}

\section{Introduction}

Recent interests have grown on connecting gene expression profiles to
survival patients' times, see e.g.~\cite{mcldata,vantdeveer}, where the aim
is to assess the influence of gene expressions on the survival
outcomes. The statistical analysis of such data faces two sorts of
problems. First, the covariates are high-dimensional: the number of
covariates is much larger than the number of observations. Second, the
survival outcomes suffers from censoring, truncation, etc. The need of
proper statistical methods to analyze such data, in particular
high-dimensional right-censored data, led in the past years to
numerous theoretical and computational contributions.

When the survival times suffer from right-censoring, the problem can
be presented as follows. For an individual $i \in \{1, \ldots, n\}$,
let $T_i$ be the time of interest (e.g. the patient survival time),
let $C_i$ be the censoring time and $X_i$ be the vector of covariates
in $\R^d$, assumed to be independent copies of $T$, $C$ and
$X=(X^1,\ldots,X^d)$. We observe $Z_i=T_i \wedge C_i$,
$\delta_i=\ind{T_i\leq C_i}$ and $X_i$ for $i=1,\ldots,n$.

The covariates vector
$X$, where both genomic outcomes and clinical data may be recorded, is
in high dimension $d \gg n$ and influences the distribution of $T$ via
its conditional hazard rate given $X=x$, defined by
\begin{eqnarray*}
  \alpha_0(t, x) = \frac{f_{T|X}(t,x)}{1 - F_{T|X}(t,x)}
\end{eqnarray*}
for $t > 0$, where $f_{T|X}$ and $F_{T|X}$ are respectively the
conditional density and distribution functions of $T$ given $X=x$. In
the following, we assume that the conditional hazard fulfills the
Aalen additive hazards model~\cite{aalen80}:
\begin{equation*}
  \alpha_0(t, x) = \lambda_0(t) + x^\top \beta_0, \quad \forall t\geq
  0,
\end{equation*}
where $\lambda_0$ is the baseline hazard function and $\beta_0$
measures the influence of the covariates on the conditional hazard
function $ \alpha_0$. In~\cite{ma2007additive}, an additive hazards
model is fitted to investigate the influence of the expression levels
of 8810 genes on the (censored) survival times of 92 patients
suffering from Mantel-Cell Lymphoma, see~\cite{mcldata} for the
data. The Aalen additive hazards model is indeed an useful alternative
to the Cox model~\cite{Cox}, in particular in situations where the
proportional hazards assumption is violated. It can also ``be seen as
a first-order Taylor series expansion of a general intensity''
(see~\cite{MartinussenScheikeBook}, p.~103).

When the aim is then to understand the influence of $X$ on the
survival time $T$, one wants to estimate $\beta_0$ based on the
observations.  In small dimension $d \ll n$ and from the data $(Z_i,
\delta_i, X_i)_{i=1,\ldots,n}$, the least-squares estimator $\hat
\beta$ of the unknown $\beta_0$ is the minimizer of the quadratic
functional
\begin{align*}
  R_n(\beta) &=\beta^\top \H_n\beta - 2 \beta^\top {\bs h}_n,
\end{align*}
where $\H_n$ is the $d \times d$ symetrical positive semidefinite
matrix with entries
\begin{equation*}
  (\H)_{j, k} = \frac 1n \sum_{i=1}^n \int_0^{Z_i} \Big(X^j_i -
  \frac{\sum_{l=1}^n X^j_l \ind{Z_l \geq t}}{\sum_{l=1}^n
    \ind{Z_l \geq t}} \Big) \Big(X^k_i -\frac{\sum_{l=1}^n X^k_l
    \ind{Z_l \geq t}}{\sum_{l=1}^n \ind{Z_l \geq t}} \Big) dt,
\end{equation*}
and where ${\bs h}_n \in \R^d$ has coordinates 
\begin{equation*}
  (\bs h_n)_j = \frac 1n \sum_{i=1}^n \delta_i \Big(X^j_i -
  \frac{\sum_{k=1}^n X^j_k \ind{Z_k \geq Z_i}}{\sum_{k=1}^n
    \ind{Z_k \geq Z_i}} \Big).
\end{equation*}
When $d \leq n$ and if $\H_n$ is full
rank, we can write
\begin{equation*}
  \hat \beta = (\H_n)^{-1} {\bs h}_n,
\end{equation*}
see also~\cite{linying} or~\cite{martinussenscheike09b}. The estimator
$ \hat \beta$ is $\sqrt{n}$-consistent and asymptotically Gaussian,
see e.g. \cite{ABGK}.

When $X$ contains genomic outcomes, one typically has $d \gg n$, and
the matrix $\H_n$ is no longer of full rank. A sparsity assumption is
then natural in this setting: we expect only a few genes to have an
influence on the survival times, so we expect $\beta_0$ to be sparse,
which means that it has only a few non-zero coordinates. Several
papers use sparsity inducing penalization in the context of survival
analysis, mainly for the Cox multiplicative risks model or the Aalen
additive risks model, we refer to~\cite{witten2010survival} for a
review. Most procedures are based on $\ell_1$-penalization, where one
considers
\begin{equation}
  \label{eq:estequ}
  \hat \beta \in \argmin_{\beta \in \R^d} \Big\{ R_n(\beta) + \lambda
  \sum_{j=1}^d w_j |\beta_j| \Big\}.
\end{equation}
The smoothing parameter $\lambda > 0$ makes the balance between
goodness-of-fit and sparsity, and the $w_j \geq 0$, $j = 1, \ldots, d$
are weights allowing for a precise tuning of the penalization. The
Lasso penalization corresponds to the simple choice $w_j = 1$, while
in the adaptive Lasso~\cite{MR2279469} one chooses $w_j = |\tilde
\beta_j|^{-\gamma}$ where $\tilde \beta_j$ is a preliminary estimator
and $\gamma > 0$ a constant. The idea behind this is to correct the
bias of the Lasso in terms of variable selection accuracy, see
\cite{MR2279469} and~\cite{zhang10} for regression analysis. The
weights $w_j$ can also be used to scale each variable at the same
level, which is suitable when some variable has a strong variance
compared to the others. As a by-product of the theoretical analysis
given in this paper, we introduce a new way of scaling the variables
using data-driven weights $\hat w_j$ in the $\ell_1$ penalization,
see~\eqref{eq:weights} below.

In the Cox proportional hazards model, $R_n(\beta)$ is the partial
likelihood (see e.g. \cite{Cox} or \cite{ABGK}), for which the Lasso,
adaptive Lasso, smooth clipped absolute deviation penalizations and
the Dantzig selector are considered, respectively,
in~\cite{tibshirani97}, \cite{MR2409726}, \cite{zhang2007adaptive},
\cite{MR1892656} and \cite{MR2779635}.

For the additive risks, \cite{ma2006additive} considers principal
component regression, \cite{ma2007additive} considers a Lasso with a
least-squares criterion that differs from the one considered here,
\cite{MR2395831,martinussenscheike09b} considers the ridge, Lasso and
adaptive Lasso penalizations and~\cite{martinussenscheike09a}
considers the partial least-squares and ridge regression
estimators. 

A serious advantage, from the computational point of view, in using
additive risks over multiplicative risks has to be
highlighted. Indeed, for the additive risks, the estimating
Equation~\eqref{eq:estequ} has a least-squares form, so that one can
apply in this case the fast Lars algorithm~\cite{MR2060166} in order
to obtain the whole path of solutions of the Lasso, as explained
in~\cite{MR2395831} for instance.  This point is particularly relevant
in practice, since one typically uses splitting techniques, such as
cross-validation, to select the smoothing parameter, or ensemble
feature methods, such as stability
selection~\cite{MeinshausenBuhlmann}, to select covariates. The
motivations and main contributions of this work are enumerated in the
following.

\emph{First motivation.}  Among the papers that propose some
mathematical analysis of the statistical properties of estimators of
the form~\eqref{eq:estequ} (upper bounds, support recovery, etc.), the
results are asymptotic in the number of observations. This can be a
problem since, in practice, one can not, in general, consider that the
asymptotic regime has been reached: in~\cite{mcldata}, for example,
the expression levels of 8810 genes and survival information
are measured for only 92 patients. Considering only the references
that are the closest to the work proposed here, the oracle property
for the adaptive Lasso is given in~\cite{MR2395831}, which is an
asymptotic property about the support and the asymptotic distribution
of the estimator, and asymptotic normality and consistency in variable
selection for the adaptive Lasso is proved
in~\cite{martinussenscheike09b}, where results about the Dantzig
selector are also derived using the restricted isometry property and
the uniform uncertainty principle from~\cite{MR2382644}. While
non-asymptotic results, like sparse oracle inequalities for instance,
are now well-known for regression or density estimation (see for
instance~\cite{MR2676897}, \cite{MR2533469},
\cite{bertin11:_adapt_dantiz}, among many others), such results are
not yet available for survival data. In this paper, we establish the
first results of this kind for survival analysis.

\emph{Second motivation.} We give sharp oracle inequalities (with
leading constant 1) for the prediction error associated to the
survival problem. The results are stated for general counting
processes, including the censoring case, while most papers consider
censored data only. Our results are stated without the assumption that
the intensity is linear in the covariates. In fact, our Lasso
estimator can be computed using an arbitrary dictionary of functions,
so that one can expect a better approximation of the true underlying
intensity.

\emph{Third motivation.} In order to prove our results, we need a new
version of Bernstein's inequality for martingales with jumps, where
the predictable variation, which is not observable in this problem, is
replaced by the optional variation, which is observable. This
concentration inequality is of independent interest, and could be
useful for other statistal problems as well. 

\emph{Fourth motivation.} Finally, and more importantly, our
non-asymptotic analysis leads to an adaptive data-driven weighting of
the $\ell_1$-norm, that involves the optional variation of each
element of the dictionary (or of each covariate in the linear
case). More precisely, our sharp control of the noise term exhibits
the fact that the $\ell_1$-penalization (see~\eqref{eq:estequ}) should
be scaled using data-driven weights of order (writing only the
dominating terms, see Section~\ref{sec:construction} for details)
\begin{equation*}
  \hat w_j \approx \sqrt{\frac{x + \log d} n \hat V_j},
\end{equation*}
where
\begin{equation*}
  \hat V_j = \frac 1n \sum_{i=1}^n \delta_i \Big(X^j_i -
  \frac{\sum_{k=1}^n X^j_k \ind{Z_k \geq Z_i}}{\sum_{k=1}^n
    \ind{Z_k \geq Z_i}} \Big)^2
\end{equation*}
corresponds, roughly, to an estimate of the variance of variable $j$.
Hence, our theoretical analysis exhibits a new way of tuning the
$\ell_1$ penalization, by multiplying each coordinate by this
empirical variance term, in order to make less apparent eventual
differences between the variability of each $X^j$ for $j = 1, \ldots,
d$. This particular form of weighting, or scaling of the variables,
was not previsouly noticed in literature.
 
The paper is organized as follows. Section~\ref{sec:model} describes
the model. The Lasso estimator is constructed
in~Section~\ref{sec:construction}. Oracle inequalities for the Lasso
are given in Section~\ref{sec:oracle-inequalities}, see
Theorems~\ref{thm:slow-oracle} and~\ref{thm:fast-oracle1}. Some
details about the construction of the least-squares criterion are
given in Section~\ref{sec:decomprisk}. The data-driven Bernstein's inequality
is stated in Section~\ref{sec:martingales}, see
Theorem~\ref{thm:empirical-bernstein}, and the proofs of our results
are given in Section~\ref{sec:proofs}.


\section{High dimensional Aalen model}
\label{sec:model}

Let $(\Omega, \mathcal F, \P)$ be a probability space and $(\mathcal
F_t)_{t \geq 0}$ a filtration satisfying the usual conditions:
increasing, right-continuous and complete
(see~\cite{JacodShiryaev}). Let $N$ be a marked counting process with
compensator $\Lambda$ with respect to $(\mathcal F_t)_{t \geq 0}$, so
that $M = N - \Lambda$ is a $(\mathcal F_t)_{t \geq 0}$-martingale. We
assume that $N$ is a marked point process satisfying the \emph{Aalen
  multiplicative intensity model}. This means that $\Lambda$ writes
\begin{equation}
  \label{eq:aalen}
  \Lambda(t) = \int_0^t \alpha_0(s, X) Y_s ds
\end{equation}
for all $t \geq 0$, where:
\begin{itemize}
\item the intensity $\alpha_0$ is an unknown deterministic and
  nonnegative function called \emph{intensity}
\item $X \in \mathbb R^d$ is a $\mathcal F_0$-measurable random vector
  called \emph{covariates} or \emph{marks};
\item $Y$ is a predictable random process in $[0, 1]$.
\end{itemize}
With differential notations, this model can be written has
\begin{equation}
  \label{eq:model}
  d N_t = \alpha_0(t, X) Y_t \,dt + d M_t
\end{equation}
for all $t \geq 0$ with the same notations as before, and taking $N_0
= 0$. Now, assume that we observe $n$ i.i.d. copies
\begin{equation}
  \label{eq:whole_sample}
  D_n = \{ (X_i, N^i_t, Y^i_t) : t \in [0, \tau], 1 \leq i \leq n \}
\end{equation}
of $\{ (X, N_t, Y_t) : t \in [0, \tau] \}$, where $\tau$ is the
end-point of the study. Without loss of generality, we set $\tau = 1$.
We can write
\begin{equation*}
  dN^i_t = \alpha_0(t, X_i) Y^i_t dt + d M^i_t
\end{equation*}
for any $i = 1, \ldots, n$ where $M^i$ are independent $(\mathcal
F_t)_{t \geq 0}$-martingales. In this setting, the random variable
$N^i_t$ is the number of observed failures during the time interval
$[0, t]$ of the individual~$i$.  This model encompasses several
particular examples: censored data, marked Poisson processes and
Markov processes, see e.g.~\cite{ABGK} for a precise exposition. In the censored case, described in the Introduction, the random processes in $D_n$, see Equation~\eqref{eq:whole_sample}, are given by
\begin{eqnarray*}
  N^i (t) = \ind{Z_i \leq t, \delta_i=1} \; \text{ and } \; Y^i (t) =
  \ind{Z_i \geq t}
\end{eqnarray*}
for $i=1,\dots,n$ and $0 \leq t \leq 1$.

In this paper, we assume that the intensity function satisfies the
Aalen additive model in the sense that it writes
\begin{equation}
  \label{eq:additive-model}
  \alpha_0(t, x) = \lambda_0(t) + h_0(x),
  \end{equation}
  where $\lambda_0 : \mathbb R_+ \to \mathbb R_+$ is a nonparametric
  \emph{baseline} intensity and $h_0 : \mathbb R^d \to \mathbb
  R_+$. Note that in the ``usual'' Aalen additive model,
  see~\cite{linying, mckeaguesasieni, martinussenscheike09a,
    martinussenscheike09b}, the function $h_0$ is linear:
  \begin{equation*}
    h_0(x)=x^\top  \beta_0,
  \end{equation*}
  where $\beta_0$ is an unknown vector in $\mathbb R^d$.  The aim of
  the paper is to recover the function $h_0$ based on the observation
  of the sample $D_n$.

\section{Construction of an $\ell_1$-penalization procedure}
\label{sec:construction}

\subsection{A least-squares type functional}

The problem considered here is a regression problem: we want to
explain the influence of the covariates $X_i$ on the survival data
$N^i$ and $Y^i$. Namely, we want to infer on $h_0$, while the baseline
function $\lambda_0$ is considered as a nuisance parameter. Thanks to
the additive structure~\eqref{eq:additive-model}, we can construct an
estimator of $h_0$ without any estimation of $\lambda_0$, so that the
influence of the covariates on the survival data can be infered
without any knowledge on $\lambda_0$. This classical principle leads
to the construction of the partial likelihood in the Cox model
(multiplicative risks, see~\cite{Cox}) and to the construction of the
``partial'' least-squares (in reference to the partial likelihood) for
the additive risks, see \cite{linying}, which is the one considered
here. The ``partial least-squares'' criterion for a ``covariate''
function $h : \R^d \rightarrow \R^+$ is defined as:
\begin{equation}
  \label{eq:risk1}
  h \mapsto \frac 1n \sum_{i=1}^n \int_0^1 ( h(X_i) - \bar h_{Y}(t)
  )^2 Y^i_t dt - \frac 2n \sum_{i=1}^n \int_0^1 (h(X_i) - \bar
  h_{Y}(t)) d N^i_t,
\end{equation}
where 
\begin{equation*}
  \bar h_Y(t) = \frac{\sum_{i=1}^n h(X_i) Y^i_t}{\sum_{i=1}^n
    Y^i_t}.
\end{equation*} 
It has been first introduced in~\cite{linying}. The main steps leading
to~\eqref{eq:risk1} are described in Section~\ref{sec:decomprisk}
below, where we explain why it is indeed suitable for the estimation
of $h_0$ (see in particular Equation~\eqref{eq:suited}).

Now, we consider a set
\begin{equation*}
  \cH = \{ h_1, \ldots, h_M \}
\end{equation*}
of functions $h_j : \R^M \rightarrow \R^+$, called \emph{dictionary},
where $M$ is large ($M \gg n$). The set $\cH$ can be a collection of
basis functions, that can approximate the unknown $h$, like wavelets,
splines, kernels, etc. They can be also estimators computed using an
independent training sample, like several estimators computed using
different tuning parameters, leading to the so-called aggregation
problem, see~\cite{MR2351101} for instance. Implicitely, it is assumed
that the unknown $h_0$ is well-approximated by a linear combination
\begin{equation}
  \label{eq:h-linear}
  h_\beta(x) = \sum_{i=1}^M \beta_j h_j(x),
\end{equation}
where $\beta \in \R^M$ is to be estimated. However, note that we won't
assume, for the statements of our results, that the unknown $h_0$ is
equal to $h_{\beta_0}$ for some unknown $\beta_0 \in \R^M$, hence
allowing for a model bias. Note that the setting considered here
includes the linear case: if $h_j(x) = x_j$ with $d = M$, then the
estimator has the form $\hat h(x) = x^\top \hat \beta$. Introducing
\begin{equation}
  \bar h_{j, Y}(t) = \frac{\sum_{i=1}^n h_j(X_i) Y^i_t}{\sum_{i=1}^n
    Y^i_t} \; \text{ and } \; \bar h_{\beta, Y}(t) = \sum_{j=1}^M
  \beta_j \bar h_{j, Y}(t),
\end{equation}
we define the least-squares risk of $\beta \in \R^M$ as
\begin{equation}
  \label{eq:risk}
  R_n(\beta) = \frac 1n \sum_{i=1}^n \int_0^1 ( h_\beta(X_i) -
  \bar h_{\beta, Y}(t) )^2 Y^i_t dt - \frac 2n \sum_{i=1}^n \int_0^1
  (h_\beta(X_i) - \bar h_{\beta, Y}(t)) d N^i_t,
\end{equation}    
which is equal to the functional~\eqref{eq:risk1} where we
applied~\eqref{eq:h-linear}. Note that~\eqref{eq:risk} is a
least-squares criterion, since
\begin{equation}
  \label{eq:risk-simple1}
  R_n(\beta) = \beta^\top \H_n \beta - 2 \beta^\top \bs h_n,
\end{equation}
where $\H_n$ is the $M \times M$ matrix with entries
\begin{equation}
  \label{eq:gram}
  (\H)_{j, k} = \frac 1n \sum_{i=1}^n \int_0^1 (h_j(X_i) - \bar
  h_{j, Y}(t) ) (h_k(X_i) - \bar h_{k, Y}(t) ) Y^i_t dt,
\end{equation}
and where $\bs h_n \in \R^M$ has coordinates
\begin{equation*}
  (\bs h_n)_j = \frac 1n \sum_{i=1}^n \int_0^1 (h_j(X_i) - \bar
  h_{j, Y}(t) )  d N_t^i.
\end{equation*}
Since $\H_n$ is a symetrical positive semidefinite matrix, we can take
\begin{equation*}
  \G_n = \H_n^{1/2},
\end{equation*}
so that
\begin{equation*}
    \label{eq:risk-simple2}
    R_n(\beta) = | \G_n \beta |_2^2 - 2 \beta^\top \bs h_n,
\end{equation*}
where $|x|_2$ stands for the $\ell_2$-norm of $x \in \R^n$. Note that
we will denote by $|x|_p$ the $\ell_p$ norm of $x$.

\subsection{$\ell_1$-penalization for the Aalen model}

For the problem considered here, we have seen that the empirical risk
$R_n$ has to be chosen with care. This is also the case for the
$\ell_1$ penalization to be used for this problem. Namely, for a
well-chosen sequence of positive data-driven weights $\hat w = (\hat
w_1, \ldots, \hat w_M)$, we consider the weighted $\ell_1$-norm
\begin{equation}
  \label{eq:pen}
  \pen(b) = |b|_{1, \hat w} = \sum_{j=1}^M \hat w_j |b_j|,
\end{equation}
and choose $\hat \beta$ according to the following penalized
criterion:
\begin{equation}
  \label{eq:lasso-def}
  \hat \beta_n \in \argmin_{b \in B} \Big\{ R_n(b) + \pen(b) \Big\}
\end{equation}
where $B$ is an arbitrary convex set (typically $B = \R^M$ or $B =
\R^M_+$, the latter making $h_{\hat \beta_n}$ non-negative). The
weights considered in~\eqref{eq:lasso-def} are given by $\hat w_j =
\hat w(h_j)$ (where we recall that $h_j \in \cH$) and where for any
function $h$, we take
\begin{equation}
  \label{eq:weights}
  \hat w(h) = c_{1} \sqrt{\frac{x + \log M + \hat \ell_{n, x}(h)} n
    \hat V(h)} + c_{2} \frac {x + 1 + \log M + \hat \ell_{n, x}(h)}
  n \norm{h}_{n, \infty},
\end{equation}
where:
\begin{itemize}
\item $x > 0$ and $c_1 = 2 \sqrt 2$, $c_2 = 4 \sqrt{14/3} + 2/3$,
\item $\norm{h}_{n, \infty} = \max_{i=1, \ldots, n} |h(X_i)|$,
\item $\hat V(h)$ is a term corresponding to the ``observable
  empirical variance'' of $h$ (see below for details), given by
  \begin{equation*}
    \hat V(h) = \frac 1n \sum_{i=1}^n \int_0^1 (h(X_i) - \bar
    h_{Y}(t))^2 d N_t^i,
  \end{equation*}
\item $\hat \ell_{n, x}(h)$ is a small technical term coming out of
  our analysis:
  \begin{equation*}
    \hat \ell_{n, x}(h) = 2 \log \log \Big( \frac{6 e n \hat V(h) + 56
      x \norm{h}_{n, \infty}^2}{24 x \norm{h}_{n, \infty}^2 } \vee e
    \Big).
  \end{equation*} 
\end{itemize}
Note that the weights $\hat w_j$ are fully data-driven. The shape of
these weights comes from a new empirical Bernstein's inequality
involving the optional variation of the noise process of the model,
see Theorem~\ref{thm:empirical-bernstein} in
Section~\ref{sec:martingales} below.

The penalization~\eqref{eq:pen} is tuned for the estimation problem at
hand. It uses the estimator $\hat V(h)$ of the (unobservable)
predictable quadratic variation
\begin{equation*}
  V(h) = \frac 1n \sum_{i=1}^n \int_0^1 (h(X_i) - \bar h_Y(t))^2
  \alpha_0(t, X_i) Y_t^i dt,
\end{equation*}
and it does not depend on an uniform upper bound for $h$. As a
consequence, it can give, from a practical point of view, some insight
into the tuning of the $\ell_1$-penalization. In particular, our
analysis prove that the $j$-th coordinate of $\beta$ in the $\ell_1$
penalization should be rescaled by $\hat V(h_j)^{1/2}$. Note that this
was not previously noticed in literature, in part because most results
are stated using an asymptotic point of view, see the references
mentioned in Introduction.

\section{Oracle inequalities}
\label{sec:oracle-inequalities}


If $\beta \in \R^M$, we denote its support by $J(\beta) = \{ j \in \{
1, \ldots, M \} : \beta_j \neq 0 \}$ and its \emph{sparsity} is
$|\beta|_0 = |J(\beta)| = \sum_{j=1}^M \ind{\beta_j \neq 0}$, where
$\ind{A}$ is the indicator of $A$ and $|B|$ is the cardinality of a
finite set $B$. If $J \subset \{ 1, \ldots, M \}$, we also introduce
the vector $\beta_J$ such that $(\beta_J)_j = \beta_j$ if $j \in J$
and $(\beta_J)_j = 0$ if $j \in J^\complement$, where $J^\complement =
\{ 1, \ldots, M \} - J$. We define the empirical norm of a function
$h$ by
\begin{equation}\label{eq:empnorm}
  \norm{h}_n^2 = \frac 1n \sum_{i=1}^n \int_0^1 ( h(X_i) - \bar
  h_{Y}(t) )^2 Y^i_t \,dt,
\end{equation}
and remark that $\norm{h_\beta}_n^2 = |\G_n \beta|_2^2 / n$.

Below are two oracle inequalities for $h_{\hat \beta}$. The first one
(Theorem~\ref{thm:slow-oracle}) is a ``slow'' oracle inequality, with
a rate of order $(\log M / n)^{1/2}$, which holds without any
assymption on the Gram matrix $\G_n$. The second one
(Theorem~\ref{thm:fast-oracle1}) is an oracle inequality with a fast
rate of order $\log M / n$, that holds under an assumption on the
restricted eigenvalues of $\G_n$. 

\begin{theorem}
  \label{thm:slow-oracle}
  Let $x > 0$ be fixed, and let $\hat h = h_{\hat \beta}$, where
  \begin{equation*}
    \hat \beta_n \in \argmin_{b \in B} \Big\{ R_n(b) + \pen(b) \Big\},
  \end{equation*}  
  with $\pen(b)$ given by~\eqref{eq:pen}. Then we have, with a
  probability larger than $1 - 29 e^{-x}$\textup:
  \begin{equation*}
    \norm{\hat h - h_0}_n^2 \leq \inf_{\beta \in B} \Big(
    \norm{h_\beta - h_0}_n^2 +  2 \pen(\beta) \Big).
  \end{equation*}
\end{theorem}
Note that
\begin{align*}
  \pen(\beta) \leq |\beta|_1 \max_{j=1, \ldots, M} \Bigg[ c_{1}
  &\sqrt{\frac{x + \log M + \hat \ell_{n, x}(h_j)} n \hat V(h_j)} \\
  &+ c_{2} \frac {x + 1 + \log M + \hat \ell_{n, x}(h_j)} n
  \norm{h_j}_{n, \infty} \Bigg]
\end{align*}
for any $\beta \in \R$, so the dominant term in $\pen(\beta)$ is, up
to the slow $\log \log$ term, of order $|\beta|_1 \sqrt{\log M / n}$,
which is the expected slow rate for $\hat h$ involving the
$\ell_1$-norm (see \cite{MR2533469} for the regression model and
\cite{MR2676897, bertin11:_adapt_dantiz} for density estimation).

For the proof of oracle inequalities with a fast $\log M / n$ rate,
the \emph{restricted eigenvalue} (RE) condition introduced
in~\cite{MR2533469} and~\cite{MR2555200,MR2500227} is of importance.
Restricted eigenvalue conditions are implied by, and in general weaker
than, the so-called \emph{incoherence} or RIP assumptions, which
excludes strong correlations between covariates. This condition is
acknowledged to be one of the weakest to derive fast rates for the
Lasso. One can find in~\cite{MR2576316} an exhaustive survey and
comparison of the assumptions used to prove fast oracle inequalities
for the Lasso, where the so-called ``compatibility condition'', which
is slightly more general than RE, is described.

The restricted eigenvalue condition is defined below. Note that our
presentation (and arguments used in the proof of
Theorem~\ref{thm:fast-oracle1}) is close
to~\cite{koltchinskii2010nuclear}, where oracle inequalities for the
matrix Lasso are given. Let us introduce, for any $\beta \in \R^M$ and
$c_0 > 0$, the cone
\begin{equation}
  \label{eq:cone}
  \C_{\beta, c_0} = \big\{ b \in \R^M : |b_{J(\beta)^\cp} |_{1,
    \hat w} \leq c_0 |b_{J(\beta)}|_{1, \hat w} \big\}.
\end{equation}
The cone $\C_{\beta, c_0}$ consists of vectors that have a support
close to the support of $\beta$. Then, introduce
\begin{equation}
  \label{eq:mu_c0}
  \mu_{c_0}(\beta) = \inf \Big\{ \mu > 0 : |b_{J(\beta)}|_2 \leq
  \frac{\mu}{\sqrt n} |\G_n b|_2 \quad \forall b \in \C_{\beta, c_0}
  \Big\}.
\end{equation}
The number $1 / \mu_{c_0}(\beta)$ is an uniform lower bound for $|\G_n
b|_2 / |b_{J(\beta)}|_2$ over $b \in \C_{\beta, c_0}$. Hence, it is a
lower bound for ``eigenvalues'' restricted over vectors with a support
close to the support of $\beta$. Also, note that $c \mapsto
\mu_c(\beta)$ is non-increasing.
\begin{theorem}
  \label{thm:fast-oracle1}
  Let $x > 0$ be fixed and let $\hat h = h_{\hat \beta}$, where
  \begin{equation*}
    \hat \beta_n \in \argmin_{b \in B} \Big\{ R_n(b) + 2 \pen(b)
    \Big\},
  \end{equation*}  
  with $\pen(b)$ given by~\eqref{eq:pen}. Then we have, with a
  probability larger than $1 - 29 e^{-x}$\textup:
  \begin{equation*}
    \norm{h_{\hat \beta} - h_0}_n^2 \leq \inf_{\beta \in B} \Big(
    \norm{h_{\beta} - h_0}_n^2 + \frac{9}{4} \mu_3(\beta)^2
    |\hat w_{J(\beta)}|_2^2 \Big),   
  \end{equation*}
  where
  \begin{equation*}
    |\hat w_{J(\beta)}|_2^2 = \sum_{j \in J(\beta)} \hat w_j^2.
  \end{equation*}
\end{theorem}
Note that
\begin{align*}
  |\hat w_{J(\beta)}|_2^2 \leq 2 |\beta|_0 \max_{j \in J(\beta)}
  \Bigg[ &c_{1}^2 \frac{x + \log M + \hat \ell_{n, x}(h_j)}{n} \hat
  V(h_j)  \\
  &+ c_{2}^2 \Big( \frac{x + 1 + \log M + \hat \ell_{n, x}(h_j)}{n}
  \norm{h_j}_{n, \infty} \Big)^2 \Bigg],
\end{align*}
so the dominant term is (up to the $\log \log$ term) of order
$|\beta|_0 \log M / n$. This is the fast rate to be found in sparse
oracle inequalities~\cite{MR2533469,MR2555200,MR2382644}. Moreover,
note that the (sparse) oracle inequality in
Theorem~\ref{thm:fast-oracle1} is sharp, in the sense that there is a
constant $1$ in front of the oracle term $\inf_{\beta \in B}
\norm{h_{\beta} - h_0}_n^2$, see Remark~\ref{rem:sharp} below.

Now, let us state Theorem~\ref{thm:fast-oracle1} under the restricted
eigenvalue condition.
\begin{assumption}[$\re(s, c_0)$ \cite{MR2533469}]
  For some integer $s \in \{ 1, \ldots, M \}$ and a constant $c_0 >
  0$, we assume that $\G_n$ satisfies:
  \begin{equation*}
  0<  \kappa(s, c_0) = \min_{\substack{J \subset \{ 1, \ldots, M \}, \\
        |J| \leq s}} \min_{\substack{b \in \R^M \setminus \{ 0 \}, \\
        |b_{J^\cp}|_{1, \hat w} \leq c_0 |b_J|_{1, \hat w} }}
    \frac{|\G_n b|_2}{\sqrt n |b_J|_2}
  \end{equation*}  
\end{assumption}

Note that using the previous notations, we have
\begin{equation*}
  \kappa(s, c_0) = \min_{\substack{b \in \R^M \setminus \{ 0 \} \\ |b|_0 \leq
      s}} \frac{1}{\mu_{c_0}(b)}.
\end{equation*}

\begin{corollary}
 \label{cor:fast-oracle2}
 Let $x > 0$, $s \in \{ 1, \ldots, M \}$ be fixed and let $\hat h$ be
 the same as in Theorem~\ref{thm:fast-oracle1}. Then, under Assumption
 $\RE(s, 3)$, we have, with a probability larger than $1 - 29
 e^{-x}$\textup:
 \begin{equation*}
   \norm{h_{\hat \beta} - h_0}_n^2 \leq \inf_{\substack{\beta \in B \\
       |\beta|_0 \leq s}} \Big( \norm{h_{\beta} - h_0}_n^2 +
   \frac{9}{4 \kappa(s, 3)^2} |\hat w_{J(\beta)}|_2^2 \Big).
 \end{equation*} 
\end{corollary}


\begin{remark}
  \label{rem:c0}
  Note that the constant $c_0 = 3$ (for $\mu_{c_0}(\beta)$) is used in
  Theorem~\ref{thm:fast-oracle1}. This is because with a large
  probability, $\hat \beta - \beta$ belongs to the cone $\C_{\beta,
    3}$. Such an argument of cone constraint is at the core of the
  convex analysis underlying the proof of fast oracle inequalities for
  the Lasso, see for instance
  \cite{MR2382644,MR2533469,koltchinskii2010nuclear}.
\end{remark}

\begin{remark}
  \label{rem:sharp}
  We were able to prove a sharp sparse oracle inequality (with leading
  constant $1$), because we adapted in our context a recent argument
  from~\cite{koltchinskii2010nuclear}, that uses some tools from
  convex analysis (such as the fact that the subdifferential mapping
  is monotone, see~\cite{MR0274683}) in the study of $\hat \beta$ as
  the minimum of the convex functional $R_n + \pen$.
\end{remark}

\section{An empirical Bernstein's inequality}
\label{sec:martingales}

The proofs of Theorems~\ref{thm:slow-oracle}
and~\ref{thm:fast-oracle1} require a sharp control of the ``noise
term'' arising from model~\eqref{eq:model}. For a fixed function $h$,
this noise term is the stochastic process
\begin{equation*}
  Z_t(h) = \frac 1n \sum_{i=1}^n \int_0^t (h(X_i) - \bar h_Y(s)) d
  M_s^i,
\end{equation*}
where we recall that $M_t^i = N_t^i - \Lambda_t^i$ are
i.i.d. martingales with jumps with jumps of size $+1$, as we assume
the existence of the intensity function $\alpha_0$,
see~\eqref{eq:aalen}. In order to give an upper bound on $|Z_t|$ that
holds with a large probability, one can use Bernstein's inequality for
martingales with jumps, see \cite{lipstershiryayev}, and note that a
proof of this fact is implicit in the proof of
Theorem~\ref{thm:empirical-bernstein}, see Section~\ref{sec:proofs}
below. Applied to the process $Z_t(h)$, this writes
\begin{equation*}
  \P\Big[ |Z_t(h)| \geq \sqrt{\frac{2 v x}{n}} + \frac{x}{3n}, V_t(h)
  \leq v \Big] \leq 2 e^{-x}
\end{equation*}
for any $x, v > 0$, where
\begin{equation*}
  V_t(h) = n \inr{Z(h)}_t = \frac 1n \sum_{i=1}^n \int_0^t (h(X_i) -
  \bar h_Y(s))^2 \alpha_0(s, X_i) Y^i_s \,ds
\end{equation*}
is the predictable variation of $Z_t$, which will also be referred to
as variance term. But, since the term $V_t(h)$ depends explicitly on
the unknown intensity $\alpha_0$, one cannot use it in the penalizing
term of the Lasso estimator. Morever, this result is stated on the
event $\{ V_t \leq v \}$ while we would like an inequality that holds
in general. Hence, we need a new Bernstein's type inequality, that
uses an observable empirical variance term instead of $V_t(h)$. We
prove in Theorem~\ref{thm:empirical-bernstein} below that we can
replace $V_t(h)$ by the optional variation of $Z_t(h)$, which can be
also seen as an estimator of $V_t(h)$ and is defined as:
\begin{equation*}
  \hat V_t(h) = n [Z(h)]_t = \frac 1n \sum_{i=1}^n
  \int_0^t (h(X_i) - \bar h_Y(s))^2 d N_s^i.
\end{equation*}
Moreover, our result holds in general, and not on $\{ V_t(h) \leq v
\}$. The counterpart for this is the presence of a small $\log \log$
term in the upper bound for $|Z_t(h)|$.
\begin{theorem}
  \label{thm:empirical-bernstein}
  For any numerical constants $c_\ell > 1, \eps > 0$ and $c_0 > 0$
  such that $e c_0 > 2(4/3 + \eps) c_\ell$, the following holds for
  any $x > 0$:
  \begin{align}
    \label{eq:empirical-bernstein}
    \P\Bigg[ |Z_t(h)| \geq c_1 \sqrt{\frac{x + \hat \ell_{n, x}(h)} n
      \hat V_t(h) } + c_2 \frac {x + 1 + \hat \ell_{n, x}(h)} n
    \norm{h}_{n, \infty} \Bigg] \leq c_3 e^{-x},
  \end{align}
  where
  \begin{equation*}
    \hat \ell_{n, x}(h) = c_\ell \log \log \Big( \frac{2 e n \hat V_t(h)
      + 8 e (4/3 + \eps) x \norm{h}_{n, \infty}^2}{4 (e c_0 - 2 (4/3 +
      \eps) c_\ell) \norm{h}_{n, \infty}^2 } \vee
    e \Big), \quad \norm{h}_{n, \infty} = \max_{i=1, \ldots, n} |h(X_i)| 
  \end{equation*}
  and where
  \begin{align*}
    c_1 &= 2 \sqrt{1 + \eps}, \quad c_2 =
    2 \sqrt{2 \max(c_0, 2 (1 + \eps) (4/3 + \eps))} + 2/3, \\
    c_3 &= 8 + 6(\log(1 + \eps))^{-c_\ell} \sum_{j \geq 1}
    j^{-c_\ell}.
  \end{align*}
  By choosing $c_\ell = 2$, $\epsilon = 1$ and $c_0 = 4 (4/3 + \eps)
  c_\ell / e = 56 / (3e)$, Inequality~\eqref{eq:empirical-bernstein}
  holds with the following numerical values:
  \begin{align*}
    c_{1} &= 2 \sqrt 2, \quad c_2 = 4 \sqrt{14 / 3} + 2/3 \leq 9.31 \\
    c_3 &= 8 + (\log 2)^{-2} \pi^2 + 4 \leq 28.55, \\
    \hat \ell_{n, x}(h) &= 2 \log \log \Big( \frac{2 e n \hat V_t(h) +
      56 e x \norm{h}_{n, \infty}^2 / 3}{8 \norm{h}_{n, \infty}^2}
    \vee e \Big).
  \end{align*}
\end{theorem}

The concentration inequality~\eqref{eq:empirical-bernstein} is fully
data-driven, since the random variable that upper bounds $|Z_t(h)|$
with a large probability is observable. Note that the numerical values
given in Theorem~\ref{thm:empirical-bernstein} are the one used in the
construction of the $\ell_1$-penalization~\eqref{eq:pen}. These are
chosen for the sake of simplicity, but another combination of
numerical values can be considered as well. 

The idea of using Bernstein's deviation inequality with an estimated
variance is of importance for statistical
problems. In~\cite{bertin11:_adapt_dantiz} for instance, a Bernstein's
inequality with empirical variance is derived in order to study the
Dantzig selector for density estimation. Note that, however, we are
not aware of a previous result such as
Theorem~\ref{thm:empirical-bernstein} for continuous time martingales
with jumps, excepted for a work in progress~\cite{wip-pat-vinc}, which
uses a similar concentration inequalities in the context of point
processes.

\subsection*{Acknowledgements}

We wish to thank the authors of \cite{wip-pat-vinc} for giving us a
preliminary version of their work, that helped in simplifying the
proof of Theorem~\ref{thm:empirical-bernstein}. We thank also Sylvain
Delattre and Sarah Lemler for helpful discussions.




\section{Proofs}
\label{sec:proofs}

\subsection{Decomposition of the least-squares}
\label{sec:decomprisk}

In this section, we give the details of the construction of the
partial least-squares~\eqref{eq:risk1}. It is based on the
decomposition, using the additive structure~\eqref{eq:additive-model},
of the least-squares risk for counting processes depending on
covariates, see for instance~\cite{Patricia2} and~\cite{CGG}. In
model~\eqref{eq:model}, on the basis of the
observations~\eqref{eq:whole_sample}, the least-squares functional to
be considered for the estimation of $\alpha_0$ is given by
\begin{equation*}
  L_n(\alpha) =  \frac 1n \sum_{i=1}^n \int_0^1 \alpha^2(t, X_i)
  Y^i_t \,dt - \frac 2n  \sum_{i=1}^n \int_0^1 \alpha(t, X_i) dN^i_t,
\end{equation*}
where $\alpha : \mathbb R_+ \times \mathbb R^d \to \mathbb R_+$. Now,
if $\alpha(t,x) = \lambda(t) + h(x)$, we can decompose $L_n$ in the
following way:
\begin{equation}
  \label{eq:decomp}
  L_n(\alpha) = L_{n, 1}(\lambda) + L_{n, 2}(h) + L_{n, 3}(\lambda, h),
\end{equation}
where
\begin{align*}  
  L_{n, 1}(\lambda) &= \frac {1}{n} \sum_{i=1}^n \int_0^1
  (\lambda(t) + \bar h_{Y}(t) )^2 Y^i_t \,dt - \frac 2n \sum_{i=1}^n
  \int_0^1 (\lambda(t) +\bar h_{Y}(t)) d N^i_t \\
  L_{n, 2}(h) &= \frac {1}{n} \sum_{i=1}^n \int_0^1 ( h(X_i) - \bar
  h_{Y}(t) )^2 Y^i_t \,dt - \frac 2n \sum_{i=1}^n \int_0^1 (h(X_i)
  - \bar h_{Y}(t)) d N^i_t \\
  L_{n, 3}(\lambda, h) &= \frac{2}{n} \sum_{i=1}^n \int_0^{1}(\lambda(t) + \bar
  h_{Y}(t)) (h(X_i) - \bar h_{Y}(t))Y^i_t\, dt,
\end{align*}
where, as introduced in~Section~\ref{sec:construction}:
\begin{equation*}
  \bar h_Y(t) = \frac{\sum_{i=1}^n h(X_i) Y^i_t}{\sum_{i=1}^n
    Y^i_t}.
\end{equation*} Now, the
point is that, according to Lemma~\ref{lem:orthotrick} below, the term
$L_{n, 3}$ is zero.
\begin{lemma}
  \label{lem:orthotrick}
  For any function $h : \mathbb R^d \rightarrow \R^+$ and any function
  $\varphi : \R^+ \rightarrow \R^+$, we have
  \begin{equation*}
    \sum_{i=1}^n \int_0^1 \varphi(t) (h(X_i) - \bar h_Y(t))  Y^i_t \,
    dt = 0.
  \end{equation*}
\end{lemma}
Lemma~\ref{lem:orthotrick} follows from an easy computation which is
omitted. The term $L_{n, 2}$ in~\eqref{eq:decomp} is the partial
least-squares criterion considered in Section~\ref{sec:construction},
see Equation~\eqref{eq:risk1}. We now explain why it is suitable for
the estimation of $h_0$. If the Aalen additive model holds, we have
$dN^i_t = (\lambda_0(t) + h_0(X_i)) Y^i_t\, dt + dM^i_t$ for all
$i=1,\ldots,n$, so we can write, using again
Lemma~\ref{lem:orthotrick}:
\begin{align*}
  L_{n,2}(h) &=\frac {1}{n} \sum_{i=1}^n \int_0^1 ( h(X_i) - \bar
  h_{Y}(t) )^2 Y^i_t\, dt \\
  & \quad - \frac 2n \sum_{i=1}^n \int_0^1 (h(X_i) - \bar h_{Y}(t))
  (h_0(X_i)- \bar h_{0,Y}(t))Y^i_t \,dt \\
  & \quad - \frac 2n \sum_{i=1}^n \int_0^1 (h(X_i) - \bar h_{Y}(t))
  dM^i_t,
\end{align*}
where 
\begin{equation*}
  \bar h_{0,Y}(t) = \frac{\sum_{i=1}^n h_0(X_i) Y^i_t}{\sum_{i=1}^n
    Y^i_t}.
\end{equation*}
Now, using the empirical norm $\norm{\cdot}_n^2$ defined in
Equation~\eqref{eq:empnorm}, see Section~\ref{sec:construction} above,
we can finally write
\begin{align}
  \label{eq:suited}
  L_{n,2}(h) = \norm{h-h_0}_n^2 - \norm{h_0}_n^2 - \frac 2n
  \sum_{i=1}^n \int_0^1 (h(X_i) - \bar h_{Y}(t)) dM^i_t.
\end{align}
The last term in the right hand side of~\eqref{eq:suited} is a noise
term, with tails controlled in Section~\ref{sec:martingales} above. It
is now understood that finding a minimizer of $L_{n,2}$, or a
penalized version of it, is a natural way of estimating $h_0$. We
refer the reader to~\cite{martinussenscheike09b} for an other
justification of the ``partial least-squares'' criterion in the linear
case $h_0(x) = x^\top \beta_0$.


\subsection{Proof of Theorem~\ref{thm:empirical-bernstein}}

For $i=1,\ldots,n$, the processes $N^i$ are i.i.d. counting processes
satisfying the Doob-Meyer decomposition $N^i_t - \int_0^t \alpha_0(s,
X_i) Y^i_s \,ds = M^i_t$, see Equation~\eqref{eq:model}. This implies
that the processes $M^i$ are i.i.d. centered martingales, with
predictable variation $\inr{M^i}_t = \int_0^t
\alpha_0(s,X_i)Y^i_s\,dt$ and optional variation $[M^i]_t = N^i_t$,
see e.g.~\cite{ABGK} for details. Moreover, the jumps of each $M^i$,
denoted by $\Delta M^i_t = M^i_t - M^i_{t_-}$, are in $\{ 0, 1\}$.
Introduce the process
\begin{equation*}
  U_t = \frac 1n \sum_{i=1}^n \int_0^t H_s^i d M_s^i
\end{equation*}
where
\begin{equation*}
  H^i_t = \frac{h(X_i) - \bar h_Y(t)}{2 \max_{i=1, \ldots, n} |h(X_i)|}.
\end{equation*}
Note that $|H_t^i| \leq 1$. Since $H^i$ is predictable and bounded,
the process $U$ is a square integrable martingale, as a sum of square
integrable martingales. Its predictable variation $\inr{U}$ is given
by:
\begin{equation*}
  \vartheta_t = n \inr{U}_t = \frac 1n \sum_{i=1}^n \int_0^t  (H_s^i)^2 
\alpha_0(s, X_i) Y^i_s \,ds,
\end{equation*} while its optional variation $ [U]$ is given by
\begin{equation*}
  \hat \vartheta_t = n [U]_t = \frac 1n \sum_{i=1}^n \int_0^t
  (H_s^i)^2 d N_s^i.
\end{equation*}
From~\cite{vandegeer95}, we know that
\begin{equation}\label{eq:supermar}
  \exp( \lambda U_t - S_\lambda(t)) 
\end{equation}
is a supermartingale if $S_\lambda$ is the compensator of
\begin{equation*}
  E_t=  \sum_{0 \leq s \leq t} \big\{ \exp( \lambda \Delta U_s) - 1 - \lambda
  \Delta U_s  \big\}.
\end{equation*}
We now derive the expression of $S_{\lambda}$. The process $E$ can
also be written as
\begin{align*}
  E_t &= \sum_{s \leq t} \sum_{k\geq 2} \frac{\lambda^k}{k!} (\Delta
  U(s))^k=\sum_{s \leq t} \sum_{k\geq 2}
  \frac{\lambda^k}{k!n^k}\Big(\Delta \big(\sum_{i=1}^n \int_0^s H_u^i
  d M_u^i\big) \Big)^k \\
  &= \sum_{s \leq t} \sum_{k\geq 2} \frac{\lambda^k}{k!n^k}
  \sum_{i=1}^n \Big(\Delta \int_0^s H_u^i d M_u^i \Big)^k,
\end{align*}
where the last inequality holds almost surely, since the $M^i$ are
independent, hence do not jump at the same time (with probability
$1$).  Now, note that
\begin{equation*} 
  \Big(\Delta \int_0^s H_u^i d M_u^i \Big)^k = (H_s^i)^k \Delta M^i(s)
  =  (H_s^i)^k \Delta N^i(s),
\end{equation*}
so that we have
\begin{equation*}
  S_\lambda(t) = \sum_{i=1}^n \int_0^t \phi \Big( \frac \lambda n
  H_s^i \Big) \alpha_0(s, X_i) Y^i_s\, ds
\end{equation*}
with $\phi(x) = e^x - x - 1$. The fact that~\eqref{eq:supermar} is a
supermartingale entails
\begin{equation}
  \label{eq:devia1}
  \P\Big[ U_t \geq \frac{S_\lambda(t)}{\lambda} + \frac{x}{\lambda}
  \Big] \leq e^{-x}
\end{equation}
for any $\lambda, x > 0$. The following facts hold true:
\begin{itemize}
\item $\phi(x h) \leq h^2 \phi(x)$ for any $0 \leq h \leq 1$ and $x >
  0$ ;
\item $\phi(\lambda) \leq \frac{\lambda^2}{2(1 - \lambda /
    3)}$ for any $\lambda \in (0, 3)$ ;
\item $\min_{\lambda \in (0, 1/b)} \big( \frac{a \lambda}{1 - b
    \lambda} + \frac x \lambda \big) = 2 \sqrt{ax} + b x,$ for any $a,
  b, x > 0$.
\end{itemize}
For any $w > 0$, they entail the following embeddings:
\begin{align}
  \nonumber
  \Big\{ U_t \geq \sqrt{\frac{2 w x}{n}} + \frac{x}{3n}, \vt_t \leq w
  \Big\} &= \Big\{ U_t \geq \frac{\lambda_w}{2(n - \lambda_w / 3)} w +
  \frac{x}{\lambda_w}, \vt_t \leq w \Big\} \\
  \nonumber
  &\subset \Big\{ U_t \geq \frac{\phi(\lambda_w / n)}{\lambda_w} n
  \vt_t + \frac{x}{\lambda_w}, \vt_t \leq w \Big\} \\
  \label{eq:embedd1}
  &\subset \Big\{ U_t \geq \frac{S_{\lambda_w}(t)}{\lambda_w} +
  \frac{x}{\lambda_w}, \vt_t \leq w \Big\},
\end{align}
where $\lambda_w$ achieves the infimum. This leads to the standard
Bernstein's inequality:
\begin{equation*}
  \P \Big[ U_t \geq \sqrt{\frac{2 w x}{n}} 
  + \frac{x}{3n}, \vartheta_t \leq w \Big] \leq e^{-x}.
\end{equation*}
By choosing $w = c_0 (x+1) / n$ for some constant $c_0 > 0$, this
gives the following inequality, which says that when the variance term
$\vt_t$ is small, the sub-exponential term is dominating in
Bernstein's inequality:
\begin{equation}
  \label{eq:deviaUtcase1}
  \P \Big[ U_t \geq \Big(\sqrt{2c_0} + \frac 13\Big) \frac {x+1}n,
  \vartheta_t \leq \frac{c_0 (x+1)}{n} \Big] \leq e^{-x}.
\end{equation}
For any $0 < v < w < +\infty$, we have
\begin{align*}
  \Big\{ U_t \geq \sqrt{\frac{2 w \vartheta_t x}{v n}} + \frac{x}{3n}
  \Big\} \cap \{ v < \vartheta_t \leq w \} \subset \Big\{ U_t \geq
  \sqrt{\frac{2 w x}{n}} + \frac{x}{3n} \Big\} \cap \{ v < \vartheta_t
  \leq w \},
\end{align*}
so, together with~\eqref{eq:devia1} and~\eqref{eq:embedd1}, we obtain
\begin{equation}
  \label{eq:devia2_Ut}
  \P \Big[ U_t \geq \sqrt{\frac{2 w \vartheta_t x}{v n}} +
  \frac{x}{3n}, v < \vartheta_t \leq w \Big] \leq e^{-x}.
\end{equation}  
Now, we want to replace $\vt_t$ by the observable $\hat \vt_t$ in the
deviation~\eqref{eq:devia2_Ut}. Note that the process $\tilde U_t$
given by
\begin{align*}
  \tilde U_t = \hat \vt_t - \vt_t &= \frac 1n \sum_{i=1}^n \int_0^t
  (H_s^i)^2 \big(dN^i_s- \alpha_0(s, X_i) Y^i_s\, ds\big) \\
  &=\frac 1n \sum_{i=1}^n \int_0^t (H_s^i)^2 dM^i_s
\end{align*} 
is a martingale, so following the same steps as for $U_t$, we obtain
that $\exp(\tilde U_t - \lambda \tilde S_\lambda(t))$ is a
supermartingale, with
\begin{equation*}
  \tilde S_\lambda(t) = \sum_{i=1}^n \int_0^t \phi \Big( \frac \lambda n
  (H_s^i)^2 \Big) \alpha_0(s, X_i) Y^i_s\, ds.
\end{equation*}
Now, writing again~\eqref{eq:embedd1} for $\tilde U_t$ with the fact
that $|H_t^i| \leq 1$ and using the same arguments as before, we
arrive at
\begin{equation*}
  \P \Big[ |\hat \vartheta_t - \vartheta_t| \geq \frac{\phi(\lambda /
    n)}{\lambda} n \vartheta_t + \frac x \lambda  \Big] \leq 2
  e^{-x}
\end{equation*}
and
\begin{equation}
  \label{eq:devia2_vvt}
  \P \Big[ |\hat \vt_t - \vt_t | \geq \sqrt{\frac{2 w \vartheta_t x}{v
      n}} + \frac{x}{3n}, v < \vartheta_t \leq w \Big] \leq 2
  e^{-x}.
\end{equation}
But, if $\vt_t$ satisfies
\begin{equation*}
  |\hat \vt_t - \vt_t| \leq \sqrt{\frac{2 w \vartheta_t x}{v n}} +
  \frac{x}{3n},
\end{equation*}
then it satisfies
\begin{equation*}
  \vt_t \leq  2 \hat \vt_t + 2\Big( \frac wv + \frac 13 \Big) \frac xn
\end{equation*}
and 
\begin{equation*}
  \hat \vt_t \leq 2 \vt_t + 2 \Big( \frac 13 + \sqrt{\frac
    w v \Big( \frac 13 + \frac wv \Big)} + \frac{2w}{v} \Big) \frac{x}{n},
\end{equation*}
simply by using the fact that $A \leq b + \sqrt{a A}$ entails $A \leq
a + 2b $ for any $a, A, b > 0$. This proves that
\begin{equation}
  \label{eq:bernstein-inclusion-trick1}  
  \begin{split}
    \Big\{ U_t \leq & \sqrt{\frac{2 w \vartheta_t x}{v n}} +
    \frac{x}{3n} \Big\} \cap \Big\{ |\hat \vt_t - \vt_t| \leq
    \sqrt{\frac{2 w \vartheta_t x}{v n}} + \frac{x}{3n} \Big\} \\
    &\subset \Big\{ U_t \leq 2 \sqrt{ \frac{w x}{v n} \hat \vt_t } +
    \Big( 2 \sqrt{\frac wv \Big(\frac wv + \frac 13 \Big) } + \frac 13
    \Big) \frac x n \Big\},
  \end{split}
\end{equation}
so using~\eqref{eq:devia2_Ut} and~\eqref{eq:devia2_vvt}, we obtain
\begin{equation*}
  \P \Big[ U_t \geq 2 \sqrt{ \frac{w x}{v n} \hat \vt_t } + \Big( 2
  \sqrt{\frac wv \Big(\frac wv + \frac 13 \Big) } + \frac 13 \Big)
  \frac x n, v \leq \vt_t < w  \Big] \leq 3 e^{-x}.
\end{equation*}
This inequality is similar to~\eqref{eq:devia2_Ut}, where we replaced
$\vt_t$ by the observable $\hat \vt_t$ in the sub-Gaussian term.  It
remains to remove the event $\{v \leq \vt_t < w\}$ from this
inequality. First, recall that~\eqref{eq:deviaUtcase1} holds, so we
can work on the event $\{ \vt_t > c_0 (x+1) / n \}$ from now on. We
use a peeling argument: define, for $j \geq 0$:
\begin{equation*}
  v_j = c_0 \frac{x+1}{n} (1 + \epsilon)^j,
\end{equation*}
and use the following decomposition into disjoint sets:
\begin{equation*}
  \{ \vt_t > v_0 \} = \bigcup_{j \geq 0} \{ v_j < \vt_t \leq v_{j+1} \}.
\end{equation*}
We have
\begin{equation*}
  \P \Big[ U_t \geq c_{1,\eps} \sqrt{ \frac{x}{n} \hat \vt_t } +
  c_{2, \eps} \frac x n, v_j < \vt_t \leq v_{j+1}  \Big] \leq 3 e^{-x},
\end{equation*}
where we introduced the constants
\begin{equation*}
  c_{1, \eps} = 2 \sqrt{1 + \eps} \text{ and } c_{2, \eps} = 2
  \sqrt{(1 + \eps) (4/3 + \eps)} + 1/3.
\end{equation*}
Let us introduce
\begin{equation*}
  \ell = c_\ell \log \log \Big(\frac{\vt_t}{v_0} \vee e\Big),
\end{equation*}
where $c_\ell > 1$. On the event
\begin{equation*}
  \Big\{ |\hat \vt_t - \vt_t| \leq
  \sqrt{\frac{2 (1+\eps) \vartheta_t (x + \ell)}{n}} + \frac{x +
    \ell}{3n} \Big\}
\end{equation*}
we have
\begin{equation*}
  \vt_t \leq 2 \hat \vt_t + 2(4/3 + \eps) \frac xn +
  \frac{2(4/3 + \eps) c_\ell}{n} \log \log( \frac{\vt_t}{v_0} \vee e),
\end{equation*}
which entails, assuming that $e c_0 > 2(4/3 + \eps) c_\ell$:
\begin{equation*}
  \vt_t \leq \frac{e c_0}{e c_0 - 2(4/3 + \eps) c_\ell} \Big(2 \hat
  \vt_t + 2 (4/3 + \eps) \frac xn \Big),
\end{equation*}
where we used the fact that $\log \log(x) \leq x / e - 1$ for any $x
\geq e$. This entails, together
with~\eqref{eq:bernstein-inclusion-trick1}, the following embedding:
\begin{align*}
  \Big\{ U_t \leq & \sqrt{\frac{2 (1+\eps) \vartheta_t (x + \ell)}{n}}
  + \frac{x + \ell}{3n} \Big\} \cap \Big\{ |\hat \vt_t -
  \vt_t| \leq \sqrt{\frac{2 (1+\eps) \vartheta_t (x + \ell)}{n}} +
  \frac{x + \ell}{3n} \Big\} \\
  \quad &\subset \bigg\{ U_t \leq c_{1,\eps} \sqrt{\frac{\hat
      \vartheta_t (x + \hat \ell)}{n}} + c_{2,\eps} \frac{x + \hat
    \ell}{n} \bigg\},
\end{align*}
where
\begin{equation*}
  \hat \ell = c_\ell \log \log \Big( \frac{2 e n \hat \vt_t + 2 e( 4/3
    + \eps) x}{e c_0 - 2(4/3 + \eps) c_\ell} \vee e \Big).
\end{equation*}
Now, using the previous embeddings together with~\eqref{eq:devia2_Ut}
and~\eqref{eq:devia2_vvt}, we obtain
\begin{align*}
  \P &\bigg[ U_t \geq c_{1,\eps} \sqrt{\frac{\hat \vartheta_t (x +
      \hat \ell)}{n}} + c_{2,\eps} \frac{x + \hat
    \ell}{n}, \vt_t > v_{0} \bigg] \\
  &\leq \sum_{j \geq 0} \P \Big[ U_t \geq \sqrt{\frac{2 (1+\eps)
      \vartheta_t (x + \ell)}{v n}} + \frac{x + \ell}{3n}, v_j < \vt_t
  \leq v_{j+1} \Big] \\
  &\quad + \sum_{j \geq 0} \P\Big[ |\hat \vt_t - \vt_t| \geq
  \sqrt{\frac{2 (1+\eps) \vartheta_t (x + \ell)}{n}} + \frac{x +
    \ell}{3n}, v_j < \vt_t \leq v_{j+1} \Big]  \\
  &\leq 3 \Big( e^{-x} + \sum_{j \geq 1} e^{-(x + c_\ell \log \log
    (v_j / v_0))} \Big) \\
  &= 3 \Big(1 + (\log(1 + \eps))^{-c_\ell} \sum_{j \geq 1} j^{-c_\ell}
  \Big) e^{-x}.
\end{align*}
Together with~\eqref{eq:deviaUtcase1}, this gives
\begin{equation*}
  \P \bigg[ U_t \geq c_{1,\eps} \sqrt{\frac{\hat \vartheta_t (x +
      \hat \ell)}{n}} + c_{3,\eps} \frac{x + 1 + \hat
    \ell}{n} \bigg] \leq \Big(4 + 3 (\log(1 + \eps))^{-c_\ell} \sum_{j
    \geq 1} j^{-c_\ell} \Big) e^{-x},
\end{equation*}
where $c_{3,\eps} = \sqrt{2 \max(c_0, 2(1+\eps)(4/3+\eps))} + 1/3$.
Now, it suffices to multiply both sides of the inequality
\begin{align*}
  U_t \geq c_{1,\eps} \sqrt{ \frac{x + \hat \ell}{n} \hat \vt_t } +
  c_{3, \eps} \frac{x + 1 + \hat \ell}{n}
\end{align*}
by $2 \norm{h}_{n, \infty}$ to recover the statement of
Theorem~\ref{thm:empirical-bernstein}. \hfill $\square$



\subsection{Some notations and preliminary results for the proof of
  the oracle inequalities}
\label{sec:proofs-proba-part}


Let us introduce the following notations. Let $\bs h(\cdot) =
(h_1(\cdot), \ldots, h_M(\cdot))^\top$ and $\bar {\bs h}_Y(\cdot) =
(\bar h_{1,Y}(\cdot), \ldots, \bar h_{M, Y}(\cdot))^\top$, so that
$h_\beta = \bs h^\top \beta$ and $\bar h_{\beta, Y} = \bar {\bs
  h}^\top_Y \beta $. We will use the notation $\inr{\cdot, \cdot}_n$
for the following ``empirical'' inner-product between to functions $h,
h' : \R^d \rightarrow \R^+$ (two ``covariates'' functions):
\begin{equation*}
  \inr{h, h'}_n = \frac 1n \sum_{i=1}^n \int_0^1 (h(X_i) - \bar
  h_{Y}(t) ) (h'(X_i) - \bar h'_{Y}(t) )  Y^i_t \,dt,
\end{equation*}
and the corresponding empirical norm:
\begin{equation*}
  \norm{h}_n^2 = \frac 1n \sum_{i=1}^n \int_0^1 (h(X_i) - \bar
  h_{Y}(t) )^2 Y^i_t \,dt.
\end{equation*}
Note that with these notations, we have:
\begin{equation*}
  \beta^\top \H_n \beta' = \inr{h_\beta, h_{\beta'}}_n.
\end{equation*}
To avoid any possible confusion, we will always write $\beta^\top
\beta'$ for the Euclidean inner product between two vectors $\beta$
and $\beta'$ in $\R^M$.

In view of~\eqref{eq:gram}, we can write
\begin{equation*}
  \H_n = \frac 1n \sum_{i=1}^n \int_0^1 ({\bs h}(X_i) - \bar
  {\bs h}_{Y}(t) ) (\bs h(X_i) - \bar {\bs h}_{Y}(t) )^\top Y^i_t\, dt,
\end{equation*}
and
\begin{equation*}
  \bs h_n = \frac 1n \sum_{i=1}^n \int_0^1 (\bs h(X_i) - \bar {\bs
    h}_{Y}(t) )  d N_t^i.
\end{equation*}
Now, in view of (\ref{eq:additive-model}) and~(\ref{eq:model}), the
following holds:
\begin{equation}
  \label{eq:hndecomp}
  \bs h_n = \bs h_n' + \bs Z_n,
\end{equation}
where:
\begin{align*}
  (\bs h_n')_j &= \frac 1n \sum_{i=1}^n \int_0^1 (h_j(X_i) - \bar
  h_{j, Y}(t) ) (\lambda_0(t) + h_0(X_i) ) Y^i_t\, dt, \\
  (\bs Z_n)_j &= \frac 1n \sum_{i=1}^n \int_0^1 (h_j(X_i) - \bar
  h_{j, Y}(t) ) d M_t^i.
\end{align*}
Using Lemma~\ref{lem:orthotrick} two times, we obtain:
\begin{align*}
  (\bs h_n')_j &= \frac 1n \sum_{i=1}^n \int_0^1 (h_j(X_i) - \bar
  h_{j, Y}(t)) (\lambda_0(t) + h_0(X_i)) Y^i_t dt \\
  &= \frac 1n \sum_{i=1}^n \int_0^1 (h_j(X_i) - \bar
  h_{j, Y}(t)) h_0(X_i) Y^i_t dt \\
  &= \frac 1n \sum_{i=1}^n \int_0^1 (h_j(X_i) - \bar h_{j, Y}(t))
  (h_0(X_i) - \bar h_{0, Y}(t)) Y^i_t dt,
\end{align*}
namely
\begin{equation}
  \label{eq:hprimetrick}
  (\bs h_n')_j = \inr{ h_j, h_0 }_n.
\end{equation}

\subsection{Proof of Theorem~\ref{thm:slow-oracle}}

Recall that the empirical risk $R_n$ is given
by~\eqref{eq:risk-simple1}. As a consequence of~\eqref{eq:hndecomp}
and~\eqref{eq:hprimetrick}, we obtain the following decomposition of
the empirical risk:
\begin{equation*}
  R_n(\beta) = \beta^\top \H_n \beta - 2 \beta^\top\bs h_n =
  \norm{h_\beta}_n^2 - 2 \inr{h_\beta, h_0}_n - 2 \beta^\top \bs Z_n,
\end{equation*}
so, for any $\beta \in \R^M$, the following holds:
\begin{equation*}
  R_n(\hat \beta) - R_n(\beta) = \norm{h_{\hat \beta} - h_0}_n^2 -
  \norm{h_\beta - h_0}_n^2 +  2 (\beta - \hat \beta)^\top \bs Z_n.
\end{equation*}
By definition of $\hat \beta$, we have 
\begin{equation*}
  R_n(\hat \beta) + \pen(\hat \beta) \leq R_n(\beta) +
  \pen(\beta)
\end{equation*}
for any $\beta \in \R^M$, so:
\begin{equation*}
  \norm{h_{\hat \beta} - h_0}_n^2 \leq \norm{h_\beta - h_0}_n^2 +  2
  (\hat \beta - \beta)^\top \bs Z_n + \pen(\beta) - \pen(\hat \beta).
\end{equation*}
Let us introduce the event
\begin{equation}
  \label{eq:eventA}
  A = \bigcap_{j=1}^M \Big\{ 2 | (\bs Z_n)_j | \leq \hat w_j \Big\},
\end{equation}
where the weights $\hat w_j$ are given by~\eqref{eq:weights}. Using
Theorem~\ref{thm:empirical-bernstein} together with an union bound, we
have that
\begin{equation*}
  \P(A) \geq 1 - c_3 e^{-x},
\end{equation*}
where $c_3$ is a purely numerical positive constant from
Theorem~\ref{thm:empirical-bernstein}. On $A$, we have
\begin{equation*}
  | 2 (\hat \beta - \beta)^\top \bs Z_n | \leq \sum_{j=1}^M \hat w_j |
  \hat \beta_j - \beta_j |= |\hat \beta -
  \beta|_{1, \hat w}, 
\end{equation*}
so recalling that $\pen(\beta) = \sum_{j=1}^M \hat w_j |\beta_j|$, we
obtain 
\begin{equation*}
  \norm{h_{\hat \beta} - h_0}_n^2 \leq \norm{h_\beta - h_0}_n^2 +  2
  \pen(\beta)
\end{equation*}
for any $\beta \in \R^M$, which is the statement of
Theorem~\ref{thm:slow-oracle}. \hfill $\square$

\subsection{Proof of Theorem~\ref{thm:fast-oracle1}}

Recall the following notation: for any $J \subset \{1, \ldots, M\}$
and $x \in \R^M$, we define the vector $x_J \in \R^M$ with coordinates
by $(x_J)_j = x_j$ when $j \in J$ and $(x_J)_j = 0$ if $j \in J^\cp$,
where $J^\cp = \{ 1, \ldots, M \} -J$. Recall that
\begin{equation}
  \label{eq:lasso-def2}
  \hat \beta \in \argmin_{b \in B} \Big\{ R_n(b) + 2 \pen(b) \Big\},
\end{equation}
where $B$ is a convex set. This proof uses arguments from
\cite{koltchinskii2010nuclear}.  We denote by $\partial \phi$ the
subdifferential mapping of a convex function $\phi$. The function $b
\mapsto R_n(b)$ is differentiable, so the subdifferential of
$R_n(\cdot) + 2 \pen(\cdot)$ at a point $b \in \R^M$ is given by
\begin{equation*}
  \partial (R_n + 2\pen)(b)  = \{ \grad R_n(b) \} + 2 \partial \pen(b)
  = \{ 2 \H_n b - 2 \bs h_n \} + 2 \partial \pen(b).
\end{equation*}
So, Equation~\eqref{eq:lasso-def2} means that there is $\hat
\beta_\partial \in
\partial \pen(\hat \beta)$ such that $\grad R_n(\hat \beta) + 2 \hat
\beta_\partial$ belongs to the normal cone of $B$ at $\hat \beta$:
\begin{equation}
  \label{eq:normal-cone}
  (2 \H_n \hat \beta - 2 \bs h_n + 2 \hat \beta_\partial )^\top (\hat
  \beta - \beta) \leq 0 \quad \forall \beta \in B.
\end{equation}
Inequality~\eqref{eq:normal-cone} can be written,
using~\eqref{eq:hndecomp} and~\eqref{eq:hprimetrick}, in the following
way:
\begin{equation*}
  2 \inr{h_{\hat \beta} - h_\beta, h_{\hat \beta} - h_0}_n + 2
  (\hat \beta_\partial - \beta_\partial)^\top (\hat \beta - \beta)
  \leq - 2 \beta_\partial^\top (\hat \beta - \beta) + 2 \bs Z_n^\top
  (\hat \beta - \beta),
\end{equation*}
where chose any $\beta_\partial \in \partial \pen(\beta)$.  Now, we
use the fact that the subdifferential mapping is monotone (this is an
immediate consequence of its definition, see \cite{MR0274683},
Chapter~24, p.~240) to say that $(\hat \beta_\partial -
\beta_\partial)^\top(\hat \beta - \beta) \geq 0$. Moreover, it is
standard to see that
\begin{equation*}
  \partial |b|_{1, \hat w} = \Big\{ e + f : e_j = \hat w_j
  \sgn(b_j) \text{ and } f_{J(b)} = 0, |f_j| \leq \hat w_j \text{ for any }
  j = 1, \ldots, M \Big\},
\end{equation*}
where $J(b) = \{ j : b_j \neq 0 \}$. Let $\beta \in B$ be fixed, and
denote $J = J(\beta) = \{ j : \beta_j \neq 0 \}$. Consider $e$ and $f$
such that $\beta_\delta = e + f$, with $e_{J^\cp} = 0$. We have
$|e^\top(\hat \beta - \beta)| \leq |\hat \beta_J - \beta_J|_{1, \hat w}$
and we can take $f$ such that $f^\top(\hat \beta - \beta) = f^\top\hat
\beta_{J^\cp} = |\hat \beta_{J^\cp}|_{1, \hat w}$. This gives
\begin{equation*}
  2 \inr{h_{\hat \beta} - h_\beta, h_{\hat \beta} - h_0}_n + 2 |\hat
  \beta_{J^\cp}|_{1, \hat w} \leq 2 |\hat \beta_J - \beta_J|_{1, \hat
    w}  + 2 \bs Z_n^\top (\hat \beta - \beta).
\end{equation*}
Using Pythagora's Theorem, we have
\begin{equation}
  \label{eq:prodca-norm}
  2 \inr{h_{\hat \beta} - h_0, h_{\hat \beta} - h_{\beta}}_n
  = \norm{h_{\hat \beta} - h_0}_n^2  + \norm{h_{\hat \beta} -
    h_\beta}_n^2 - \norm{h_\beta - h_0}_n^2,
\end{equation}
so
\begin{align*}
  \norm{h_{\hat \beta} - h_0}_n^2 + &\norm{h_{\hat \beta} -
    h_\beta}_n^2 + 2 |\hat \beta_{J^\cp}|_{1, \hat w} \\
  &\leq \norm{h_{\hat \beta} - h_0}_n^2 + 2 |\hat \beta_J -
  \beta_J|_{1, \hat w} + 2 \bs Z_n^\top (\hat \beta - \beta).
\end{align*}
If $\inr{h_{\hat \beta} - h_0, h_{\hat \beta} - h_{\beta}}_n < 0$, we
have $\norm{h_{\hat \beta} - h_0}_n < \norm{h_{\beta} - h_0}_n$, which
entails the Theorem, so we assume that $\inr{h_{\hat \beta} - h_0,
  h_{\hat \beta} - h_{\beta}}_n \geq 0$. In this case
\begin{equation*}
  2 |\hat
  \beta_{J^\cp}|_{1, \hat w} \leq 2 \inr{h_{\hat \beta} - h_0, h_{\hat
      \beta} - h_{\beta}}_n  + 2 |\hat \beta_{J^\cp}|_{1, \hat w} \leq 2
  |\hat \beta_J - \beta_J|_{1, \hat w} + 2 \bs Z_n^\top
  (\hat \beta - \beta),
\end{equation*}
which entails, together with the fact that, on $A$
(see~\eqref{eq:eventA}), we have
\begin{equation*}
  2 | \bs Z_n^\top (\hat \beta - \beta) | = 2| (\bs Z_n)_J^\top (\hat
  \beta_J - \beta_J)| + 2 |(\bs Z_n)_{J^\cp}^\top \hat \beta_{J^\cp}|
  \leq  |\hat \beta_J - \beta_J|_{1, \hat w} + |\hat
  \beta_{J^\complement}|_{1, \hat w},
\end{equation*}
that
\begin{equation*}
  |\hat \beta_{J^\complement}|_{1, \hat w} \leq 3 |\hat
  \beta_J - \beta_J|_{1, \hat w}.
\end{equation*}
This means that $\hat \beta - \beta \in \C_{\beta, 3}$
(see~\eqref{eq:cone}). So, using~\eqref{eq:mu_c0}, we have
\begin{equation}
  \label{eq:useofRE}
  |\hat \beta_J - \beta_J |_2 \leq \mu_3(\beta) | \G_n (\hat
  \beta - \beta)|_2.
\end{equation}
Note that, on $A$, we have:
\begin{equation*}
  \norm{h_{\hat \beta} - h_0}_n^2 + \norm{h_{\hat \beta} - h_\beta}_n^2
  + |\hat \beta_{J^\complement}|_{1, \hat w} \leq \norm{h_{\beta} -
    h_0}_n^2 + 3 |\hat \beta_{J} - \beta_J|_{1, \hat w}.
\end{equation*}
A consequence of~\eqref{eq:useofRE} is
\begin{equation*}
  |\hat \beta_J - \beta_J|_{1, \hat w} \leq |\hat w_J|_2 |\hat \beta_J
  - \beta_J|_2 \leq \mu_3(\beta) |\hat w_J|_2 |\G_n (\hat \beta -
  \beta) |_2,
\end{equation*}
so we arrive at
\begin{equation*}
  \norm{h_{\hat \beta} - h_0}_n^2 \leq \norm{h_{\beta} - h_0}_n^2 + 3
  \mu_3(\beta) |\hat w_J|_2 \norm{h_{\hat \beta} - h_\beta}_n  -
  \norm{h_{\hat \beta} - h_\beta}_n^2,
\end{equation*}
and finally
\begin{equation*}
  \norm{h_{\hat \beta} - h_0}_n^2 \leq \norm{h_{\beta} - h_0}_n^2 +
  \frac 94 \mu_3(\beta)^2 |\hat w_J|_2^2,
\end{equation*}
using the fact that $a x - x^2 \leq a^2 / 4$ for any $a, x >
0$. \hfill $\square$

\bibliographystyle{plain}

\bibliography{biblio}

\end{document}